wordprocessor
\documentclass[titlepage,12pt]{article}
\usepackage{amssymb}
\usepackage{latexsym}

\newcommand{\ben}{\begin{enumerate}}
\newcommand{\een}{\end{enumerate}}
\newcommand{\ble}{\begin{lem}}
\newcommand{\ele}{\end{lem}}
\newcommand{\bth}{\begin{thm}}
\renewcommand{\eth}{\end{thm}}
\newcommand{\bpr}{\begin{prop}}
\newcommand{\epr}{\end{prop}}
\newcommand{\bco}{\begin{cor}}
\newcommand{\eco}{\end{cor}}
\newcommand{\bcon}{\begin{conj}}
\newcommand{\econ}{\end{conj}}
\newcommand{\bde}{\begin{defn}}
\newcommand{\ede}{\end{defn}}
\newcommand{\bex}{\begin{exa}}
\newcommand{\eex}{\end{exa}}
\newcommand{\barr}{\begin{array}}
\newcommand{\earr}{\end{array}}
\newcommand{\btab}{\begin{tabular}}
\newcommand{\etab}{\end{tabular}}
\newcommand{\beq}{\begin{equation}}
\newcommand{\eeq}{\end{equation}}
\newcommand{\bea}{\begin{eqnarray*}}
\newcommand{\eea}{\end{eqnarray*}}
\newcommand{\bce}{\begin{center}}
\newcommand{\ece}{\end{center}}
\newcommand{\bpi}{\begin{picture}}
\newcommand{\epi}{\end{picture}}
\newcommand{\bfi}{\begin{figure} \begin{center}}
\newcommand{\efi}{\end{center} \end{figure}}

\newcommand{\bsl}{\begin{slide}{}}
\newcommand{\esl}{\end{slide}}

\newcommand{\bib}{thebibliography}

\newcommand{\qed}{\rule{1ex}{1ex}}
\newcommand{\Qed}{\rule{1ex}{1ex} \medskip}
\newcommand{\mc}[3]{\multicolumn{#1}{#2}{#3}}
\newcommand{\Mc}[1]{\multicolumn{#1}{c}{}}

\newcommand{\ol}{\overline}

\newcommand{\hso}[1]{\hspace{-1pt}}

\newcommand{\qmq}[1]{\quad\mbox{#1}\quad}

\newcommand{\emp}{\emptyset}

\newcommand{\case}[4]{\left\{\barr{ll}#1&\mbox{#2}\\#3&\mbox{#4}\earr\right.
}

\newcommand{\al}{\alpha}
\newcommand{\be}{\beta}

\newcommand{\la}{\lambda}

\renewcommand{\th}{\theta}

\newcommand{\bbS}{{\mathbb S}}

\newcommand{\cA}{{\cal A}}

\newcommand{\cP}{{\cal P}}
\newcommand{\cPD}{{\cal PD}}

\newcommand{\cS}{{\cal S}}
\newcommand{\cT}{{\cal T}}

\newcommand{\cw}{\mathop{\rm cw}}

\newcommand{\Mat}{\mathop{\rm Mat}\nolimits}

\newcommand{\tr}{\mathop{\rm tr}}

\newcommand{\wt}{\mathop{\rm wt}\nolimits}

\newcommand{\scl}{\scriptstyle}

\newcommand{\scz}{\scriptsize}

\newcommand{\oup}{Oxford University Press}

\newcommand{\U}{\mathop{U}}
\newcommand{\Z}{\mathop{Z}}
\newcommand{\C}{\mathbb{C}}
\newcommand{\ZZ}{\mathbb{Z}}

\newcommand{\gl}{\mathfrak{gl}}

\newcommand{\di}{\partial}
\newcommand{\ot}{\otimes}
\newcommand{\ts}{\,}
\newcommand{\Proof}{\noindent{\bf Proof.}\ \ }

\newtheorem{thm}{Theorem}[section]
\newtheorem{prop}[thm]{Proposition}
\newtheorem{cor}[thm]{Corollary}
\newtheorem{lem}[thm]{Lemma}
\newtheorem{conj}[thm]{Conjecture}
\newtheorem{exa}[thm]{Example}
\newtheorem{defn}[thm]{Definition}

\begin{document}
\title{A Littlewood-Richardson Rule for\\ factorial Schur functions
}
\author{Alexander I. Molev\\
Centre for Mathematics and its Applications\\
Australian National University\\
Canberra, ACT 0200\\
AUSTRALIA\\
molev@pell.anu.edu.au\\[5pt]
Bruce E. Sagan\\
Department of Mathematics\\ 
Michigan State University\\
East Lansing, MI 48824-1027\\
U.S.A.\\
sagan@math.msu.edu}

\maketitle

\pagebreak

\section{Introduction}					\label{i}

As $\lambda$ runs over all partitions with length
$l(\lambda)\leq n$,
the Schur polynomials $s_{\lambda}(x)$ form
a distinguished basis in the algebra of symmetric polynomials in
the independent variables $x=(x_1,\dots,x_n)$. By definition,
$$
s_{\lambda}(x)=
\frac{\det (x_j^{\lambda_i+n-i})_{1\leq i,j\leq n}}
{\prod_{i<j}(x_i-x_j)\qquad}.
$$
Equivalently, these polynomials can be defined by the combinatorial
formula
\beq
s_{\lambda}(x)=\sum_{T}
\prod_{\alpha\in\lambda}
x_{T(\alpha)},\label{s(x)}
\eeq
summed over semistandard tableaux $T$ of shape $\lambda$ with
entries in the set $\{1,\dots,n\}$, where $T(\alpha)$ 
is the entry of $T$ in the
cell $\alpha$. 

Any product
$s_{\lambda}(x)s_{\mu}(x)$ can be expanded as a linear
combination of Schur polynomials:
\beq							\label{LR}
s_{\lambda}(x)s_{\mu}(x)=\sum_{\nu} c_{\lambda\mu}^{\nu}\ts s_{\nu}(x).
\eeq
The classical Littlewood-Richardson rule~\cite{lr:gca}
gives a method for computing the coefficients $c_{\lambda\mu}^{\nu}$.
These same coefficients appear in the expansion of a skew Schur
function
$$
s_{\nu/\la}(x)=\sum_\mu c_{\la\mu}^\nu s_\mu(x).
$$
A number of different proofs and
variations of this rule can be found in the literature; see,
e.g.~\cite{mac:sfh,sag:sgr}, and the references therein.

To state the rule, we introduce the following notation. If $T$ is a tableau
then let $\cw(T)$ be the (reverse) {\it column word} of $T$, namely
the sequence obtained by reading the entries of $T$
from top to bottom in successive columns starting from the right-most
column.  We will call the associated total order on the cells of $T$ {\it
column order} and write $\al<\be$ if cell $\al$ comes before cell
$\be$ in this order.
A word $w=a_1\cdots a_N$ in the symbols $1,\dots,n$ is a {\it lattice
permutation\/} if for $1\leq r\leq N$ and $1\leq i< n$ the number
of occurrences of $i$ in $a_1\cdots a_r$ is at least as large as the number
of occurrences of $i+1$.

The Littlewood-Richardson rule says that the coefficient
$c_{\lambda\mu}^{\nu}$ is equal to the number of semistandard tableaux
$T$ of the shape $\nu/\mu$ and weight $\lambda$ such that $\cw(T)$ is
a lattice permutation. (One usually uses {\it row words} in the formulation
of the rule.
However, it is known that these two versions are equivalent~\cite{fg:lrm}.)
In particular, $c_{\lambda\mu}^{\nu}$ is zero unless
$\lambda,\mu\subseteq\nu$
and $|\nu|=|\lambda|+|\mu|$.

We will now state an equivalent formulation of the
Littlewood-Richardson rule~\cite{jp:sss,zel:glr,kr:bac} and establish
some notation to be used in Section~\ref{cc}.
Let $R$ denote a sequence of diagrams
$$
\mu=\rho^{(0)}\to\rho^{(1)}\to
\ldots\to\rho^{(l-1)}\to\rho^{(l)}=\nu,
$$
where $\rho\to\sigma$ means that $\rho\subset\sigma$ with
$|\sigma/\rho|=1$.
Let $r_i$ denote the row number of $\rho^{(i)}/\rho^{(i-1)}$.
Then the sequence $r_1\ldots r_l$ is called the {\it Yamanouchi
symbol} of $R$.  Equivalently, $R$ corresponds to a standard
tableau $T$ of shape $\nu/\mu$ where $r_i$ is the row number of the
entry $i$ in $T$.   A semistandard tableau $T$ {\it fits} $\nu/\mu$ if
$\cw(T)$ is the Yamanouchi symbol for some standard Young
tableau of shape $\nu/\mu$.  For example, 
$$
T=\barr{ccc} 	1&1&2\\
		2&3\earr
$$
fits $(4,3,1)/(2,1)$ since $\cw(T)=21312=r_1\ldots r_5$ corresponds to
the standard tableau
$$
\barr{cccc}	 & &2&4\\
		 &1&5\\
		3\earr
$$
or equivalently to the shape sequence
$$
R:\mu=(2,1)\to(2,2)\to(3,2)\to(3,2,1)\to(4,2,1)\to(4,3,1)=\nu.
$$
The coefficient
$c_{\lambda\mu}^{\nu}$ is then equal to the number of 
semistandard tableaux $T$ of shape $\lambda$ that fit $\nu/\mu$.

The factorial Schur function $s_{\lambda}(x|a)$ 
is a polynomial in $x$ and a doubly-infinite
sequence of variables $a=(a_i)$. 
It can be defined as the ratio of two alternants~(\ref{det}) by
analogy with the ordinary case. 
This approach goes back to Lascoux~\cite{las:pe}.
The $s_{\lambda}(x|a)$ are also a special case of the double Schubert
polynomials introduced by Lascoux and Sch\"utzenberger as explained
in~\cite{las:i}. 
The combinatorial definition~(\ref{sxa}) for the particular sequence $a$
with $a_i=i-1$
is due to Biedenharn and Louck~\cite{bl:ncs}
while the case for general $a$ is due to Macdonald~\cite{mac:sft} 
and Goulden--Greene~\cite{gg:ntr}.
The equivalence of~(\ref{det}) and~(\ref{sxa}) was established 
in~\cite{mac:sft}.

Specializing $a_i=0$ for all $i$ the functions 
$s_{\lambda}(x|a)$  turn into $s_{\lambda}(x)$. They
form a basis in the symmetric polynomials
in $x$ over $\C[a]$ so one can define the corresponding
Littlewood-Richardson coefficients $c_{\lambda\mu}^{\nu}(a)$,
see~(\ref{FLR}).  Our main result is Theorem~\ref{lrab} which
gives a combinatorial rule for calculating a two-variable
generalization $c_{\theta\mu}^{\nu}(a,b)$ of these
coefficients~(\ref{cab}), where $\theta$ is a skew diagram. 
In the case $|\nu|=|\theta|+|\mu|$ the rule turns into a rule for
computing the intertwining number of the skew
representations of the symmetric group corresponding to the diagrams
$\theta$ and $\nu/\mu$~\cite{jp:sss,zel:glr}.
Specializing further to  $\mu=\emptyset$ (respectively $\theta=\lambda$) we
get
the classical Littlewood-Richardson rule in the first (respectively
second)
formulation above.  
A completely different rule for calculating
the $c_{\la\mu}^\nu(a)$ follows from the work of Billey and
Shimozono~\cite{bs:nas} on multiplying Schubert polynomials.
A Pieri rule for multiplication of a double Schubert polynomial
by a complete or elementary symmetric polynomial is given by
Lascoux and Veigneau~\cite{vei:t}.  Lascoux has pointed out that
the Newton interpolation formula in several variables~\cite{ls:in} can
also be used to give an alternative proof of the factorial
Littlewood-Richardson rule.

In Section~\ref{mrc} we consider the specialization $a_i=i-1$.
The corresponding  coefficients $c_{\lambda\mu}^{\nu}(a)$ turn out to
be the structure constants for the center of the universal
enveloping algebra for the Lie algebra $\gl(n)$
and for an algebra of invariant differential operators in certain
distinguished bases.  We also obtain a formula which relates
these coefficients to the dimensions of skew diagrams. This implies
a symmetry property of these coefficients.

\section{Preliminaries}		\label{p}

Let $x=(x_1,\dots,x_n)$ be a finite sequence of variables and
let $a=(a_i)$, $i\in\mathbb{Z}$ be a doubly-infinite variable
sequence. The {\it generalized factorial Schur function\/}
for a partition $\lambda$ of length at most $n$ can be defined
as follows~\cite{mac:sft}.  Let
$$
(y|a)^k=(y-a_1)\cdots (y-a_k)
$$
for each $k\geq 0$. Then 
\beq				\label{det}
s_{\lambda}(x|a)=
\frac{\det \bigl[(x_j|a)^{\lambda_i+n-i}\bigr]_{1\leq i,j\leq n}}
{\prod_{i<j}(x_i-x_j)\qquad\quad}. 
\eeq
There is an explicit combinatorial formula for $s_{\lambda}(x|a)$
analogous to~(\ref{s(x)}):
\beq				\label{sxa}
s_{\lambda}(x|a)=\sum_{T}\prod_{\alpha\in\lambda}
(x_{T(\alpha)}-a_{T(\alpha)+c(\alpha)}),
\eeq
where $T$ runs over all semistandard tableaux of shape $\lambda$
with entries in $\{1,\dots,n\}$, 
$T(\alpha)$ is the entry of $T$ in the cell
$\alpha\in\lambda$ and $c(\alpha)=j-i$ is the content of
$\alpha=(i,j)$.

The highest homogeneous component of $s_{\lambda}(x|a)$ in $x$ 
obviously coincides with $s_{\lambda}(x)$. Therefore
the polynomials $s_{\lambda}(x|a)$ form a basis for
$R[x]^{S_n}$ where $R=\mathbb C[a]$, and one can define
Littlewood-Richardson type coefficients $c_{\lambda\mu}^{\nu}(a)$ by
\beq
s_{\lambda}(x|a)s_{\mu}(x|a)=
\sum_{\nu} c_{\lambda\mu}^{\nu}(a) s_{\nu}(x|a). \label{FLR}
\eeq
Comparing the highest homogeneous components in $x$ on the both sides
and using the Littlewood-Richardson Rule for the $s_\la(x)$ we see that 
\beq						\label{high}
c_{\lambda\mu}^{\nu}(a)=
\case{c_{\lambda\mu}^{\nu}}{if $|\nu|=|\lambda|+|\mu|$,}
{0}{if $|\nu|>|\lambda|+|\mu|$.}
\eeq
Contrary to the classical case, the coefficients $c_{\lambda\mu}^{\nu}(a)$
turn out to be nonzero if $|\nu|<|\lambda|+|\mu|$ and
$\lambda,\mu\subseteq\nu$. This makes it possible to compute them using
induction on $|\nu/\mu|$ while keeping $\lambda$ fixed.

The starting point of our calculation is the fact that the polynomials
$s_{\lambda}(x|a)$ possess some (characteristic) vanishing properties;
see~\cite{sah:sci,oko:qih}. We use the following result
from~\cite{oko:qih}.
For a partition $\rho$ with $l(\rho)\leq n$ define an $n$-tuple
$a_{\rho}=(a_{\rho_1+n},\dots,a_{\rho_n+1})$.

\bth[Vanishing Theorem]			\label{VT}
Given partitions $\la,\rho$ with $l(\la),l(\rho)\leq n$
$$
s_{\lambda}(a_{\rho}|a)=
\case{0}{if $\lambda\not\subseteq\rho,$}
{\prod_{(i,j)\in\lambda}(a_{\lambda_i+n-i+1}-a_{n-\lambda^t_j+j})}
{if $\la=\rho$.}
$$
where $\lambda^t$ is the diagram conjugate to $\lambda$.
\eth

In particular, $s_{\lambda}(a_{\lambda}|a)\ne 0$ for any
specialization of the sequence $a$ such that $a_i\ne a_j$ if $i\ne j$. 
We reproduce the proof of the Vanishing Theorem from~\cite{oko:qih}
for completeness. 
\bigskip

\Proof The $ij$-th entry of the determinant
in the numerator of the right hand side of~(\ref{det})
for $x=a_{\rho}$ is
\beq				\label{arho}
(a_{\rho_j+n-j+1}-a_1)\cdots
(a_{\rho_j+n-j+1}-a_{\lambda_i+n-i}).
\eeq
The condition $\lambda\not\subseteq\rho$ implies that there exists
an index $k$ such that
$\rho_k<\lambda_k$. Then for $i \leq k\leq j$ we have
$$
1\leq \rho_j+n-j+1\leq \rho_k+n-k+1\leq
\lambda_k+n-k\leq \lambda_i+n-i,
$$
and so all the entries~(\ref{arho}) with $i \leq k\leq j$ are zero. Hence,
the determinant is zero which proves the first part of the theorem.

Let us set now $x=a_{\lambda}$ in~(\ref{det}). Then the $ij$-th entry of
the determinant is
$$
(a_{\lambda_j+n-j+1}-a_1)\cdots (a_{\lambda_j+n-j+1}-a_{\lambda_i+n-i}),
$$
which equals zero for $i<j$ and is nonzero for $i=j$. This means that
the matrix is lower triangular with nonzero diagonal elements.
Taking their product and dividing by 
$$
\prod_{i<j}(a_{\lambda_i+n-i+1}-a_{\lambda_j+n-j+1})
$$
we get  the  desired equation. \hfill\qed

\section{Calculating the coefficients}			\label{cc}

We will be able to prove more general results by introducing a second
infinite sequence of variables 
denoted $b=(b_i)$, $i\in\ZZ$.  Let $\theta$ and $\mu$ be skew and
normal (i.e., skewed by $\emp$) diagrams, respectively.
Define Littlewood-Richardson type coefficients $c_{\theta\mu}^{\nu}(a,b)$
by the formula
\beq						\label{cab}
s_{\theta}(x|b)s_{\mu}(x|a)=
\sum_{\nu} c_{\theta\mu}^{\nu}(a,b) s_{\nu}(x|a),
\eeq
where $s_{\theta}(x|b)$ is defined as in~(\ref{sxa}) with $\la$
replaced by $\th$ and $a$ replaced by $b$.

As in Section~\ref{i}, consider a sequence of diagrams
\beq							\label{R}
R:\mu=\rho^{(0)}\to\rho^{(1)}\to
\ldots\to\rho^{(l-1)}\to\rho^{(l)}=\nu,
\eeq
and let $r_i$ be the row number of $\rho^{(i)}/\rho^{(i-1)}$.
Construct the set $\cT(\theta,R)$ of semistandard $\theta$-tableaux
$T$ with entries from $\{1,\ldots,n=|x|\}$ such that $T$ contains
cells $\al_1,\ldots,\al_l$ with
$$
\al_1<\ldots<\al_l\qmq{and}T(\al_i)=r_i,\ 1\le i\le l,
$$
where $<$ is column order.
Let us distinguish the entries in $\al_1,\dots,\al_l$ by barring
each of them.  For example, if $n=2$ and
$$
R:(2,1)\to(2,2)\to(3,2)
$$
so that $r_1r_2=21$ then for $\th=(3,2)/(1)$ we have
$$
\cT(\th,R)=\left\{
\barr{ccc}
	&1	&1\\
\ol{1}	&\ol{2}
\earr,
\barr{ccc}
	&\ol{1}	&\ol{2}\\
1	&2
\earr,
\barr{ccc}
	&1	&\ol{2}\\
\ol{1}	&2
\earr,
\barr{ccc}
	&1	&2\\
\ol{1}	&\ol{2}
\earr,
\barr{ccc}
	&\ol{1}	&\ol{2}\\
2	&2
\earr
\right\}.
$$
We also let
$$
\cT(\th,\nu/\mu)=\biguplus_R \cT(\th,R),
$$
where the union is over all sequences $R$ of the form~(\ref{R}).
Finally, for each cell $\alpha$ with $\al_i<\al<\al_{i+1}$,
$0\le i\le l$, set $\rho(\alpha)=\rho^{(i)}$.  (Inequalities involving
cells with out-of-range subscripts are ignored.)  For instance, if
$l=|\nu/\mu|=2$ then the following schematic diagram gives the layout of
the $\rho(\al)$
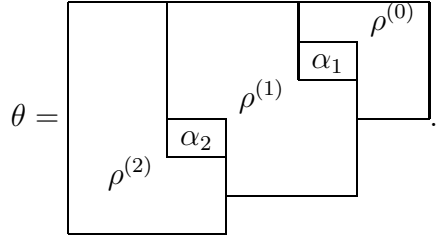
$$
\th=
\barr{|cc|c|c|c|c|} \hline
&	&\Mc{1}	&		&\Mc{1}	&\rho^{(0)}\\	\cline{5-5}
&	&\Mc{1}	&		&\al_1	&	\\	\cline{5-5}
&	&\Mc{1}	&\mc{1}{c}{\rho^{(1)}}&	&     	\\ 	\cline{3-3}\cline{6-6}
&	&\al_2	&\Mc{1}		&	\\	\cline{3-3}
&\mc{1}{c}{\rho^{(2)}}&&\Mc{1}	&	\\	\cline{4-5}
&\Mc{1}	&	\\	\cline{1-3}
\earr.
$$

We are now in a position to state the Littlewood-Richardson rule for
the $c_{\theta\mu}^{\nu}(a,b)$.
The reader should compare the following formula with the combinatorial
one for the $s_\la(x|a)$ in~(\ref{sxa}).

\bth				\label{lrab}
The coefficient $c_{\theta\mu}^{\nu}(a,b)$ is zero
unless $\mu\subseteq\nu$. If $\mu\subseteq\nu$ then
$$
c^{\nu}_{\theta\mu}(a,b)=\sum_{T\in\cT(\theta,\nu/\mu)}
\prod_{\scl\alpha\in\theta\atop \scl T(\alpha)\mbox{\scz unbarred}}
\left((a_{\rho(\alpha)})_{T(\alpha)}-b_{T(\alpha)+c(\alpha)}\right).
$$
\eth

As immediate specializations of this result, note the following.
\ben
\item If $a=b$ and $\theta$ is normal then this is a
Littlewood-Richardson rule for the $s_\la(x|a)$.
\item If $a=b$ and $\mu$ is empty then this is a rule for the expansion of
skew factorial Schur polynomial.
\item If $|\nu|=|\theta|+|\mu|$ then $c_{\theta\mu}^{\nu}(a,b)$ is
independent
of $a$ and $b$ and equals the number of semistandard tableaux of shape
$\theta$ that fit $\nu/\mu$. This coincides with the number of pictures
between
$\theta$ and $\nu/\mu$~\cite{jp:sss,zel:glr}.   In particular,
$c_{\theta\mu}^{\nu}(a,b)=c_{\theta\mu}^{\nu}$, an ordinary
Littlewood-Richardson coefficient.
\item If $\mu=\emp$ and $\theta=\la$ is normal then this is a rule for
the re-expansion of a factorial Schur polynomial in terms of those
for a different sequence of second variables. In particular,
$$
s_{\lambda}(x|a)=\sum_{\nu\subseteq\lambda}g_{\lambda\nu}(a)
s_{\nu}(x)
$$
where
$$
g_{\lambda\nu}(a)=(-1)^{|\lambda/\nu|}
\sum_{T\in\cT(\lambda,\nu)}
\prod_{\scl\alpha\in\lambda\atop \scl T(\alpha)\mbox{\scz unbarred}}
a_{T(\alpha)+c(\alpha)}.
$$
A different expression for $g_{\lambda\nu}(a)$ in terms of double
Schubert polynomials is provided by
the Newton interpolation formula in several variables~\cite{ls:in}.
\een 

We present the proof of Theorem~\ref{lrab} as a chain of propositions.

Note that the first claim of the theorem
follows immediately from the Vanishing Theorem. Indeed,
let $\nu$ be minimal (with respect to containment) among all
partitions in~(\ref{cab}) such that $c_{\theta\mu}^{\nu}(a,b)\ne 0$. 
Suppose
$\nu\not\supseteq \mu$. Then setting $x=a_{\nu}$ in~(\ref{cab}) and using
the first part of the Vanishing Theorem gives
$$
0=c_{\theta\mu}^{\nu}(a,b)s_{\nu}(a_{\nu}|a).
$$
But by the Vanishing Theorem's second part we have
$s_{\nu}(a_{\nu}|a)\ne 0$ and so a contradiction to
$c_{\theta\mu}^{\nu}(a,b)\ne 0$.

We shall assume hereafter that $\mu\subseteq\nu$ and also write
$$
|a_{\rho}|=a_{\rho_1+n}+\cdots+a_{\rho_n+1}.
$$

\bpr				\label{s/s}
If $\mu\subseteq\nu$ with $|\nu/\mu|=l$ then
$$
\frac{s_{\mu}(a_{\nu}|a)}{s_{\nu}(a_{\nu}|a)}
=\sum_{\mu\to\rho^{(1)}\to\cdots\to \rho^{(l-1)}\to\nu}
\frac{1}{(|a_{\nu}|-|a_{\rho^{(0)}}|)
\cdots (|a_{\nu}|-|a_{\rho^{(l-1)}}|)},
$$
where $\rho^{(0)}=\mu$.
\epr

\Proof Setting $x=a_{\mu}$ in~(\ref{cab}) and using the Vanishing Theorem
gives
\beq					\label{low}
c_{\theta\mu}^{\mu}(a,b)=s_{\theta}(a_{\mu}|b). 
\eeq
Further, for $\theta=(1)$ and $a=b$ relation~(\ref{cab}) takes the form
(cf.~\cite[Theorem 9.1]{oo:ssf})
$$
s_{(1)}(x|a)s_{\mu}(x|a)=s_{(1)}(a_{\mu}|a)
s_{\mu}(x|a)+\sum_{\mu\to\rho}s_{\rho}(x|a)
$$
which follows from~(\ref{low}), (\ref{high}), and the Branching Theorem
for the ordinary Schur functions.

Setting $x=a_{\nu}$ in the previous equation and using the Vanishing
Theorem we get 
\beq					\label{int}
s_{(1)}(a_{\nu}|a)s_{\mu}(a_{\nu}|a)=s_{(1)}(a_{\mu}|a)
s_{\mu}(a_{\nu}|a)+\sum_{\mu\to\rho\subseteq\nu}s_{\rho}(a_{\nu}|a). 
\eeq
We have
$$ 
s_{(1)}(a_{\nu}|a)-s_{(1)}(a_{\mu}|a)=
|a_{\nu}|-|a_{\mu}|
$$
and so~(\ref{int}) gives
$$
\frac{s_{\mu}(a_{\nu}|a)}{s_{\nu}(a_{\nu}|a)}=\frac{1}{|a_{\nu}|-|a_{\mu}|}
\sum_{\mu\to\rho\subseteq\nu}\frac{s_{\rho}(a_{\nu}|a)}{s_{\nu}(a_{\nu}|a)}.
$$
Induction on $|\nu/\mu|$ completes the proof.\hfill\Qed

It will be convenient to have a notation for sums like those
occurring in the previous proposition.  So let
\beq					\label{H}
H(\mu,\rho)=
\sum_{\mu\to\rho^{(1)}\to\cdots\to \rho^{(r-1)}\to\rho}
\frac{1}{(|a_{\rho}|-|a_{\rho^{(0)}}|)
\cdots (|a_{\rho}|-|a_{\rho^{(r-1)}}|)}, 
\eeq
and
\beq					\label{H'}
H'(\rho,\nu)=
\sum_{\rho\to\rho^{(r+1)}\to\cdots\to \rho^{(l-1)}\to\nu}
\frac{1}{(|a_{\rho}|-|a_{\rho^{(r+1)}}|)
\cdots (|a_{\rho}|-|a_{\rho^{(l)}}|)},
\eeq
where $\rho^{(0)}=\mu$ and $\rho^{(l)}=\nu$.

\bpr					\label{HH'}
We have the formula
$$
c_{\theta\mu}^{\nu}(a,b)=\sum_{\mu\subseteq\rho\subseteq\nu}
s_{\theta}(a_{\rho}|b) H(\mu,\rho) H'(\rho,\nu).
$$
\epr

\Proof We use induction on $|\nu/\mu|$, noting that~(\ref{low}) is the
base case $|\nu/\mu|=0$.
Set $x=a_{\nu}$ in~(\ref{cab}) and divide both sides by
$s_{\nu}(a_{\nu}|a)$.
By Proposition~\ref{s/s} we get
$$
c_{\theta\mu}^{\nu}(a,b)=s_{\theta}(a_{\nu}|b) H(\mu,\nu)
-\sum_{\sigma\subset\nu}c_{\theta\mu}^{\sigma}(a,b) H(\sigma,\nu).
$$
By the induction hypotheses we can write this as
\bea
c_{\theta\mu}^{\nu}(a,b)&=&s_{\theta}(a_{\nu}|b) H(\mu,\nu)
-\sum_{\sigma\subset\nu}\sum_{\mu\subseteq\rho\subseteq\sigma}
s_{\theta}(a_{\rho}|b) H(\mu,\rho)H'(\rho,\sigma)H(\sigma,\nu)\\
	&=&s_{\theta}(a_{\nu}|b) H(\mu,\nu)
-\sum_{\mu\subseteq\rho\subset\nu}s_{\theta}(a_{\rho}|b) H(\mu,\rho)
\sum_{\rho\subseteq\sigma\subset\nu}H'(\rho,\sigma)H(\sigma,\nu).
\eea
To complete the proof we note that
$$
\sum_{\rho\subseteq\sigma\subseteq\nu}H'(\rho,\sigma)H(\sigma,\nu)=0,
$$
which follows from the identity
$$
\sum_{i=1}^k\frac{1}{(u_1-u_2)\cdots(u_1-u_i)
(u_k-u_i)\cdots(u_k-u_{k-1})}=0,
$$
which holds for any variables $u_1,\dots,u_{k}$ by induction on $k>1$.
(In the denominator an empty product is, as usual, equal to 1.)\hfill\Qed

Note that a different expression for the $c_{\theta\mu}^{\nu}(a,b)$
in terms of divided differences  
can be deduced from the Newton
interpolation formula in several variables~\cite{ls:in}.

\bpr
We have the recurrence relation
\beq						\label{rr}
c_{\theta\mu}^{\nu}(a,b)=\frac{1}{|a_{\nu}|-|a_{\mu}|}\left(
\sum_{\mu\to\mu'}c_{\theta\mu'}^{\nu}(a,b)
-\sum_{\nu'\to\nu}c_{\theta\mu}^{\nu'}(a,b)\right). 
\eeq
\epr

\Proof By Proposition~\ref{HH'} it suffices to check that
$$
H(\mu,\rho) H'(\rho,\nu)=\frac{1}{|a_{\nu}|-|a_{\mu}|}\left(
\sum_{\mu\to\mu'}H(\mu',\rho) H'(\rho,\nu)
-\sum_{\nu'\to\nu}H(\mu,\rho) H'(\rho,\nu')\right).
$$
This follows from the relations
$$
\sum_{\mu\to\mu'}H(\mu',\rho)=(|a_{\rho}|-|a_{\mu}|)H(\mu,\rho)
$$
and
$$
\sum_{\nu'\to\nu}H'(\rho,\nu')=(|a_{\rho}|-|a_{\nu}|)H'(\rho,\nu).\qquad\Qed
$$

Given a sequence 
$$
R: \mu=\rho^{(0)}\to\rho^{(1)}\to
\ldots\to\rho^{(l-1)}\to\rho^{(l)}=\nu,
$$
and an index $k\in\{1,\dots,l\}$ introduce a set of $\theta$-tableaux
$\cT_k(\theta,R)$ having entries from the set $\{1,\ldots,n\}$ as follows.
Each tableau $T\in \cT_k(\theta,R)$ contains cells
$\al_1,\dots,\al_{k-1},\al_{k+1},\dots,\al_l$ such that
$$
\al_1<\ldots<\al_{k-1}<\al_{k+1}<\ldots<\al_l\qmq{and}T(\al_i)=r_i,\ 1\le
i\le l,i\neq k.
$$
As usual, we distinguish the entries in  the $\al_i, i\neq k$, by
barring them.  Now modify the $\rho(\al)$ for $R$ by defining, for 
cells with unbarred entries,
$$
\rho^{+}(\alpha)=
\case{\rho^{(k)}}{if $\al_{k-1}<\al<\al_{k+1}$,}{\rho(\al)}{otherwise,}
$$
and
$$
\rho^{-}(\alpha)=
\case{\rho^{(k-1)}}{if $\al_{k-1}<\al<\al_{k+1}$,}{\rho(\al)}{otherwise.}
$$
Also define corresponding products
\bea
\cS(R)&=&
\sum_{T\in\cT(\th,R)}\prod_{\scl\al\in\th\atop\scl T(\al)\mbox{\scz
unbarred}}
\left((a_{\rho(\alpha)})_{T(\alpha)}-b_{T(\alpha)+c(\alpha)}\right),\\
\cS_k^{+}(R)&=&
\sum_{T\in\cT_k(\th,R)}\prod_{\scl\al\in\th\atop\scl T(\al)\mbox{\scz
unbarred}}
\left((a_{\rho^{+}(\alpha)})_{T(\alpha)}-b_{T(\alpha)+c(\alpha)}\right),
\eea
and similarly for $\cS_k^-(R)$.  So Theorem~\ref{lrab} is equivalent to
\beq							\label{SR}
c_{\th\mu}^\nu(a,b)=\sum_R \cS(R).
\eeq

\bpr
Given a sequence $R$ we have
\beq							\label{S+-}
\cS(R)=\frac{1}{|a_{\nu}|-|a_{\mu}|}
\sum_{k=1}^l \left(\cS_k^{+}(R)-\cS_k^{-}(R)\right). 
\eeq
\epr

\Proof  It suffices to show that for each $k$ we have
$$
\cS_k^{+}(R)-\cS_k^{-}(R)=(|a_{\rho^{(k)}}|-|a_{\rho^{(k-1)}}|)\cS(R).
$$
Formula (\ref{S+-}) will then follow from the relation
$$
\sum_{k=1}^l(|a_{\rho^{(k)}}|-|a_{\rho^{(k-1)}}|)=|a_{\nu}|-|a_{\mu}|.
$$
For a given $T\in\cT_k(\th,R)$ the factors in the formulas for
$\cS_k^{+}(R)$
and $\cS_k^{-}(R)$ are identical except for the case where
$\al_{k-1}<\al<\al_{k+1}$ and $T(\alpha)=r_k$. 
To see what happens when we divide $\cS_k^{+}(R)-\cS_k^{-}(R)$ by
$$
|a_{\rho^{(k)}}|-|a_{\rho^{(k-1)}}|=(a_{\rho^{(k)}})_{r_k}-(a_{\rho^{(k-1)}}
)_{r_k},
$$
fix $T$ and consider its contribution to the quotient.   We need the
following easily proved formula, where we are thinking of
$u=(a_{\rho^{(k)}})_{r_k}$, $v=(a_{\rho^{(k-1)}})_{r_k}$ and
$m_i=b_{T(\al)+c(\al)}$ as $\al$ runs over 
all cells of $T$ with $\al_{k-1}<\al<\al_{k+1}$ and $T(\al)=r_k$:
$$
\frac{\prod_{i=1}^s (u-m_i)-\prod_{i=1}^s (v-m_i)}{u-v}=
\sum_{j=1}^s (u-m_1)\cdots\widehat{(u-m_j)}(v-m_{j+1})\cdots(v-m_s)
$$
(a hat indicates the factor is to be omitted). The right-hand side
of this expression can now be interpreted as the contribution to
$\cS(R)$ of all tableaux gotten from $T$ by barring one of the $r_k$
between
$\al_{k-1}$ and $\al_{k+1}$ in column order.\hfill \Qed

We now prove~(\ref{SR}) by induction on $|\nu/\mu|$.
Equation~(\ref{low}) takes care of the case $|\nu/\mu|=0$.
By the induction hypothesis,
$$
\sum_{\mu\to\mu'}c_{\theta\mu'}^{\nu}(a,b)=\sum_R \cS_1^{+}(R)
\qmq{and}
\sum_{\nu'\to\nu}c_{\theta\mu}^{\nu'}(a,b)=\sum_R \cS_l^{-}(R).
$$
So formulas~(\ref{rr}), (\ref{S+-}) and 
the following proposition complete the proof of~(\ref{SR}) and hence
Theorem~\ref{lrab}. 

\bpr
We have
$$
\sum_R \sum_{k=1}^{l-1}\cS_k^{-}(R)=\sum_R \sum_{k=2}^l\cS_k^{+}(R).
$$
\epr

\Proof We can rewrite this formula as follows:
\beq						\label{wt}
\sum_{R,k,T}\wt^{-}(R,k,T)=\sum_{R',k',T'}\wt^{+}(R',k',T'),
\eeq
where $T\in\cT_k(\theta,R)$, $k=1,\dots,l-1$ and
$T'\in\cT_{k'}(\theta,R')$, $k'=2,\dots,l$ with weights defined by
$$
\wt^-(R,k,T)=
\prod_{\alpha\in\theta\atop T(\alpha)\mbox{\scriptsize unbarred}}
\left((a_{\rho^{-}(\alpha)})_{T(\alpha)}-b_{T(\alpha)+c(\alpha)}\right)
$$
and similarly define $\wt^+(R',k',T')$.  To prove~(\ref{wt})
we will construct a bijection $(R,k,T)\longleftrightarrow (R',k',T')$
preserving the weights in the sense that
$\wt^{-}(R,k,T)=\wt^{+}(R',k',T')$.
There are three cases.
\medskip

{\sl Case 1.}
Suppose that the skew diagram $\rho^{(k+1)}/\rho^{(k-1)}$ consists of
two cells in different rows and columns. Then $R'$ is the sequence
obtained from $R$ by replacing $\rho^{(k)}$ by the other diagram 
${\rho'}^{(k)}$ such that $\rho^{(k-1)}\to {\rho'}^{(k)}\to \rho^{(k+1)}$
while $k'=k+1$ and $T'=T$.
\medskip

{\sl Case 2.}
Let $\rho^{(k+1)}/\rho^{(k-1)}$ have two cells in the same row.
Then $R'=R$, $k'=k+1$ and $T'=T$.
\medskip

{\sl Case 3.}
Let $\rho^{(k+1)}/\rho^{(k-1)}$ have two cells in the same
column, say in rows $r$ and $r+1$.  Let $(i+1,j)=(i_1,j_1)$ be the cell
of $T$ containing the correponding $\ol{r+1}$.  If there is an $r$ in
cell $(i,j)$ then it must be unbarred.  In this case let $T'$ be $T$
with the bar moved from the $r+1$ to the $r$, $R'=R$ and $k'=k+1$.  Weights
are preserved since $(a_{\rho^{(k-1)}})_r=(a_{\rho^{(k+1)}})_{r+1}$ and
$T(\al)+c(\al)$ is invariant under the change.

Now suppose cell $(i,j)$ of $T$ contains an entry  less than $r$
or $(i,j)\not\in\th$ and
let $j'=j'_1$ be the column of the left-most $r+1$ in row $i+1$. 
Since this subcase is more complicated than the others, the reader may
wish to follow along with the example given after the end of this proof.
Let $s$ be the maximum integer such that for $1\le t\le s$ we have
\ben
\item there is an $\ol{r+t}$ in cell $(i+t,j_t)$ for some $j_t$
corresponding to a 
cell in the same column as those of $\rho^{(k+1)}/\rho^{(k-1)}$,
\item if $(i+t,j'_t)$ contains the left-most $r+t$ in row $i+t$ then
$(i+t,j_t)$ is between $(i+t-1,j_{t-1})$ and $(i+t-1,j'_{t-1})$ in
column order.  (Assume this is true vacuously when $t=1$.)
\een
Note that the condition on cell $(i,j)$ implies that none of the
$r+t$'s to the left of the one in $(i+t,j_t)$ can be barred.  

We now form $T'$ by moving the bar in cell $(i+t,j_t)$ to cell $(i+t,j'_t)$
and replacing the $r+t$'s in cells
$(i+t,j'_t),(i+t,j'_t+1),\ldots,(i+t,j_t)$
by $r+t-1$'s.  Note that the result will still be a semistandard
tableau because of the assumption about $(i,j)$ and the choice of
elements to decrease.  Since the
elements from the given column of $\rho^{(k+1)}/\rho^{(k-1)}$ are
still added in the correct order in $T'$ it determines a valid $R'$,
complete except for the step where a cell is added in row $r+s$ of that
column which should be done immediately following the addition of
$r+s-1$.  Then $k'$ is the position of this $r+s$.  Invariance of
weights follows from considerations like those in the first
subcase, noting that the contribution to $\wt^-$ of each entry decreased in
$T$
is the same as that of the element on its right to $\wt^+$
in $T'$.

The inverse of this construction is similar and left to the reader.
This completes the proof of the Theorem~\ref{lrab}.\hfill\Qed

As an example of the last subcase, suppose we have the $R$ sequence
$$
(3,2,2,2)\to(3,3,2,2)\to(3,3,3,2)\to(4,3,3,2)\to(4,3,3,2,1)\to(4,3,3,3,1)
$$
with Yamanouchi symbol $r_1\ldots r_5=23154$.  Let $k=1$ so $r=2$ and
consider
$$
T=\barr{cccccc}
1	&1	&1	&1	&\ol{1}	&1\\
2	&3	&3	&3	&3	&\ol{3}\\
4	&4	&\ol{4}	&4	&\ol{5}
\earr.
$$
Then $(i+1,j)=(2,6)$ and $s=2$ with $r+1,r+2=3,4$ so
$$
T'=\barr{cccccc}
1	&1	&1	&1	&\ol{1}	&1\\
2	&\ol{2}	&2	&2	&2	&2\\
\ol{3}	&3	&3	&4	&\ol{5}
\earr.
$$
Column reading the barred elements of $T'$ and inserting $r+2=4$ after
$r+1=3$ gives the Yamanouchi symbol $15234$ corresponding to the $R'$
sequence
$$
(3,2,2,2)\!\to\!(4,2,2,2)\!\to\!(4,2,2,2,1)\!\to\!(4,3,2,2,1)\!\to\!(4,3,3,2
,1)
\!\to\!(4,3,3,3,1)
$$
and $k'=5$.

\section{Multiplication rules for Capelli operators and quantum
immanants} 			\label{mrc}

Let $E_{ij}$, $i,j=1,\dots,n$ denote the standard basis
of the general linear Lie algebra $\gl(n)$. Denote by $\Z(\gl(n))$
the center of the universal enveloping algebra $\U(\gl(n))$.
Given $\kappa=(\kappa_1,\dots,\kappa_n)\in\C^n$ consider
a $\gl(n)$-module $L(\kappa)$ of
highest weight $\kappa$. That is, $L(\kappa)$ is generated by
a nonzero vector $v$ such that
\bea
E_{ii}\cdot v&=&\kappa_i\ts v,\qmq{for}i=1,\dots,n;\\
E_{ij}\cdot v&=&0 \qmq{for}1\leq i<j\leq n.
\eea
Any element $z\in \Z(\gl(n))$ acts as a scalar 
$\omega(z)=\omega_{\kappa}(z)$ in $L(\kappa)$
and this scalar is independent of
the choice of the highest weight module $L(\kappa)$. Moreover, $\omega(z)$
is a symmetric polynomial in the shifted variables $x_1,\dots,x_n$
where $x_i=\kappa_i+n-i$. The mapping $z\mapsto\omega(z)$ defines
an algebra isomorphism
$$
\omega:\ts \Z(\gl(n))\to \C[x]^{S_n}
$$
called the {\it Harish-Chandra isomorphism\/}; 
see e.g. Dixmier~\cite[Section 7.4]{dix:ae}.

For any positive integer $m$ consider the natural action
of the complex Lie group $GL(n)$ in the algebra $\cP$
of polynomials on the vector space $\C^n\ot\C^m$.
The corresponding Lie algebra $\gl(n)$ then acts by differential operators
$$
\pi(E_{ij})=\sum_{a=1}^m x_{ia}\di_{ja},
$$
where the $x_{ia}$ are the coordinates on $\C^n\ot\C^m$ and the 
$\di_{ia}$ are the corresponding partial derivatives.
This representation is uniquely extended to an algebra
homomorphism
$$
\pi:\ts \U(\gl(n))\to\cPD
$$
where $\cPD$ is the algebra of polynomial coefficient differential
operators in the $x_{ia}$. The image of
$\Z(\gl(n))$ under $\pi$ is contained in the subalgebra $\cPD^G$
of differential operators invariant with respect to the action
of the group $G=GL(n)\times GL(m)$. Moreover, if $m\geq n$
then this restriction is an algebra isomorphism which can be called
the {\it Capelli isomorphism\/}; see~\cite{how:rci,hu:cdd} for further
details.
So if $m\geq n$ we have the triple isomorphism
$$
\C[x]^{S_n}\stackrel{\omega}{\longleftarrow}\Z(\gl(n))
\stackrel{\pi}{\longrightarrow} \cPD^G.
$$
Distinguished bases in the three algebras which correspond to each
other under
these isomorphisms were constructed in~\cite{oko:qih} 
(see also~\cite{naz:yci,oko:ybw,mol:fss}).

In the algebra $\C[x]^{S_n}$ the basis is formed by the polynomials
$s_{\lambda}(x|a)$ with $l(\lambda)\leq n$ and the sequence $a$
specialized to $a_i=i-1$ for all $i\in\ZZ$. We shall denote these
polynomials by $s^*_{\lambda}(x)$.
Explicitly~\cite{bl:ncs},
$$
s^*_{\lambda}(x)=\sum_{T}\prod_{\alpha\in\lambda}
(x_{T(\alpha)}-T(\alpha)-c(\alpha)+1), 
$$
where $T$ runs over all semistandard tableaux of shape $\lambda$
with entries in $\{1,\dots,n\}$.
We shall denote by $f_{\lambda\mu}^{\nu}$ the coefficient
$c_{\lambda\mu}^{\nu}(a)$ in this specialization of $a$. In other words,
the $f_{\lambda\mu}^{\nu}$ can be defined by the formula
$$
s^*_{\lambda}(x)s^*_{\mu}(x)=\sum_{\nu}f_{\lambda\mu}^{\nu}\ts
s^*_{\nu}(x).
$$

To describe the basis in $\Z(\gl(n))$ we introduce some more notation.
Given $k$ matrices $A$, $B$, $\dots$, $C$ of size $p\times q$
with entries from an algebra $\cA$
we regard their tensor product $A\ot B\ot\cdots\ot C$
as an element
$$
\sum A_{a_1i_1}B_{a_2i_2}\cdots C_{a_ki_k}
\ot e_{a_1i_1}\ot e_{a_2i_2}\ot\cdots\ot e_{a_ki_k}\in\cA\ot
(\Mat_{pq})^{\ot k},
$$
where $\Mat_{pq}$ denotes the space of complex $p\times q$-matrices
and the $e_{ai}$ are the standard matrix units.
The symmetric group $S_k$ acts in a natural way in the tensor space
$(\C^n)^{\ot k}$,
so that we can identify permutations from $S_k$ with
elements of the algebra $(\Mat_{nn})^{\ot k}$.

For a diagram $\lambda$ with $l(\lambda)\leq n$
denote by $T_0$ the $\lambda$-tableau obtained by filling in the cells
by the numbers $1,\dots,k=|\lambda|$ 
from left to right in successive rows starting from the
first row. We let $R_{\lambda}$ and $C_{\lambda}$ denote
the row symmetrizer and column antisymmetrizer of $T_0$ respectively.
By $c_{\lambda}(r)$ we denote the content of the cell occupied by $r$.
Introduce the matrix $E=(E_{ij})$
whose $ij$-th entry is the generator $E_{ij}$ and set
$$
\bbS_{\lambda}=\frac{1}{h(\lambda)}\tr 
(E-c_{\lambda}(1))\ot\cdots\ot (E-c_{\lambda}(k))\cdot
R_{\lambda}C_{\lambda},
$$
where the trace is taken over all the tensor factors $\Mat_{nn}$, and
$h(\lambda)$ is the product of the hook-lengths of the cells of $\lambda$:
$$
h(\lambda)=\prod_{\alpha\in\lambda}h_{\alpha}.
$$ 
The elements $\bbS_{\lambda}$ with $l(\lambda)\leq n$
form a basis in the algebra $\Z(\gl(n))$. In~\cite{oko:qih} they
were called the {\it quantum immanants\/}.

Let us now describe the basis in the algebra $\cPD^G$.
The representation $\pi$ can be written in a matrix form as follows
$$
\pi:\ts E\mapsto XD^t,
$$
where $X$ and $D$ are the $n\times m$ matrices formed
by the coordinates $x_{ia}$ and the derivatives $\di_{ia}$, respectively,
while $D^t$ is the matrix transposed to $D$.
Introduce the following differential operators
$$
\Delta_{\lambda}=\frac{1}{k!}\ts
\tr X^{\ot k}\cdot (D^t)^{\ot k}\cdot \chi^{\lambda},
$$
where $\chi^{\lambda}$ is the irreducible character of $S_k$
corresponding to $\lambda$. Explicitly,
$$
\Delta_{\lambda}=\frac{1}{k!}\sum_{i_1,\dots,i_k}\sum_{a_1,\dots,a_k}
\sum_{s\in S_k}\chi^{\lambda}(s)x_{i_1a_1}\cdots x_{i_ka_k}
\di_{i_{s(1)}a_1}\cdots \di_{i_{s(k)}a_k}.
$$
The operators $\Delta_{\lambda}$ with $l(\lambda)\leq n$ form a basis
in $\cPD^G$. They are called the {\it higher Capelli operators\/}.

The following identities were proved in~\cite{oko:qih} 
(for other proofs see~\cite{naz:yci,oko:ybw,mol:fss}):
$$ 
\omega(\bbS_{\lambda})=s^*_{\lambda}(x) \qmq{and}
\pi(\bbS_{\lambda})=\Delta_{\lambda}.
$$

Using Theorem~\ref{lrab} we obtain the following multiplication rules for
the
elements $\bbS_{\lambda}$ and the operators $\Delta_{\lambda}$.

\bth
We have
$$
\bbS_{\lambda}\bbS_{\mu}=\sum_{\nu}f_{\lambda\mu}^{\nu}\ts\bbS_{\nu}
$$
and
$$
\Delta_{\lambda}\Delta_{\mu}=\sum_{\nu}f_{\lambda\mu}^{\nu}\ts\Delta_{\nu}
$$
where the coefficients $f_{\lambda\mu}^{\nu}$ are given by
\beq						\label{f}
f_{\lambda\mu}^{\nu}=\sum_{T\in\cT(\lambda,\nu/\mu)}
\prod_{\scl\alpha\in\theta\atop\scl T(\alpha)\mbox{\scz unbarred}}
\left(\rho(\alpha)_{T(\alpha)}+n-2\ts T(\alpha)-c(\alpha)+1\right)
\eeq
with $R$, $\cT(\lambda,\nu/\mu)$, and $\rho(\alpha)$ defined in
Theorem~\ref{lrab}.  \hfill\qed
\eth

Proposition~\ref{HH'} enables us to obtain another formula for 
$f_{\lambda\mu}^{\nu}$. For a skew diagram $\nu/\mu$ let
$$
h(\nu/\mu)=\frac{|\nu/\mu|!}{\dim \nu/\mu},
$$
where $\dim \nu/\mu$ is the number of standard $\nu/\mu$-tableaux.
In particular, if $\mu$ is empty $h(\nu)$ 
is the product of the hook-lengths of the cells of $\nu$.

In the specialization of the sequence $a$ under consideration we obtain
from~(\ref{H}) and~(\ref{H'}) that
$$
H(\mu,\rho)=\frac 1{h(\rho/\mu)} 
\qmq{and}
H'(\rho,\nu)=\frac {(-1)^{|\nu/\rho|}}{h(\nu/\rho)}.
$$
Moreover, by Proposition~\ref{s/s},
$$
\frac{s_{\lambda}(a_{\rho}|a)}{s_{\rho}(a_{\rho}|a)}=
\frac 1{h(\rho/\lambda)}
$$
and by the Vanishing Theorem
$$
s_{\rho}(a_{\rho}|a)=h(\rho).
$$
Thus Proposition~\ref{HH'} becomes the following.

\bpr
One has the formula
\beq							\label{fh}
f^{\nu}_{\lambda\mu}=\sum_{\lambda,\mu\subseteq\rho\subseteq\nu}
(-1)^{|\nu/\rho|}\frac{h(\rho)}{h(\nu/\rho)h(\rho/\lambda)h(\rho/\mu)}.
\qquad\qed
\eeq
\epr

Formula~(\ref{fh}) implies the following symmetry property of the
coefficients $f^{\nu}_{\lambda\mu}$.

\bco
If $l(\lambda^t),\ l(\mu^t),\ l(\nu^t)\le n$ then
$$
f^{\nu^t}_{\lambda^t\mu^t}=f^{\nu}_{\lambda\mu}.
$$
\eco

\Proof This follows immediately from the relation
$h(\nu^t/\mu^t)=h(\nu/\mu)$.\hfill\Qed

\noindent
{\bf Remarks.} 1. It follows from~(\ref{f}) that
the coefficients $f^{\nu}_{\lambda\mu}$ are integers while the
summands on
right hand side of~(\ref{fh}) need not be.  In fact the numbers
$h(\nu/\mu)$ need not be integers either, e.g., $h((3,2)/(1))=24/5$. 

2. Note that since in the case of $|\nu|=|\lambda|+|\mu|$ the
$f^{\nu}_{\lambda\mu}$ coincide with the classical
Littlewood-Richardson coefficients $c^{\nu}_{\lambda\mu}$,
the latter can be computed using~(\ref{fh}) as well, but this does
not appear to be very useful for practical purposes. For example,
consider $\lambda=\mu=(1^n)$ and
$\nu=(2^r1^{n-r})$, then~(\ref{fh}) gives
$$
f^{\nu}_{\lambda\mu}=\sum_{k=0}^r(-1)^{r-k}\frac{(n+1)!}{k!(r-k)!(n-k+1)}
$$
while directly from~(\ref{f}) we get
\beq						\label{fact}
f^{\nu}_{\lambda\mu}=(n-r)! 
\eeq
\medskip

As a final example, take $m=n$ in the definition of $\Delta_{\lambda}$.
Then for $\lambda=(1^n)$ we get
the classical Capelli operator~\cite{cap:otf}:
$$
\Delta_{(1^n)}=\det X\det D.
$$
We find from~(\ref{f}) that the coefficients $f^{\nu}_{(1^n)(1^n)}$
are zero except for $\nu=(2^r1^{n-r})$, $r=0,1,\dots,n$ .
So by~(\ref{fact}) the square of the Capelli operator is given by
$$
(\det X\det D)^2=\sum_{r=0}^n (n-r)!\ts
\Delta_{(2^r1^{n-r})}.
$$

\medskip

{\it Acknowledgment.}
We would like to thank G.~D.~James, A.~N.~Kirillov, A.~Lascoux and S.~Sahi
for
valuable discussions. 

\begin{\bib}{99}

\bibitem[BL]{bl:ncs} L. Biedenharn and J. Louck, A new class of symmetric
polynomials defined in terms of tableaux, {\it Advances in Appl.\ Math.}
{\bf 10} (1989), 396--438.

\bibitem[BS]{bs:nas} S. Billey and M. Shimozono, A new approach to
Schubert calculus, in preparation.

\bibitem[C]{cap:otf} A. Capelli, Sur les op\'erations dans la th\'eorie des
formes alg\'ebriques,  {\it Math.\  Ann.\/} {\bf 37} (1890), 1--37.

\bibitem[D]{dix:ae} J. Dixmier, ``Alg\`ebres enveloppantes,''
Gauthier-Villars, Paris, 1974.

\bibitem[FG]{fg:lrm} S. Fomin and C. Greene, A Littlewood--Richardson
miscellany, {\it European J. Combin.\/} {\bf 14}  (1993), 191--212.

\bibitem[GG]{gg:ntr} I. Goulden and C. Greene, A new tableau representation
for supersymmetric Schur functions, {\it J. Algebra} {\bf 170} (1994),
687--703.

\bibitem[H]{how:rci} R. Howe,  Remarks on classical invariant theory,
{\it Trans.\ Amer.\  Math.\ Soc.\/} {\bf 313} (1989),  539--570.

\bibitem[HU]{hu:cdd} R. Howe and T. Umeda, The Capelli identity, the
double commutant theorem, and multiplicity-free actions, {\it Math.\ 
Ann.\/}
{\bf 290} (1991), 569--619.

\bibitem[JP]{jp:sss} G. D. James and M. H. Peel, Specht series for
skew representations of symmetric groups, {\it J. Algebra} {\bf 56} (1979),
343--364.

\bibitem[KR]{kr:bac} A. N. Kirillov and N. Yu. Reshetikhin,
The Bethe ansatz and the combinatorics of Young tableaux,
{\it J. Soviet Math.\ } {\bf 41} (1988), 925--955.

\bibitem[L1]{las:pe} A. Lascoux, 
Puissances ext\'erieures, d\'eterminants et cycles 
de Schubert, {\it Bull S.M.F.} {\bf 102} (1974), 161--179.

\bibitem[L2]{las:i} A. Lascoux, ``Interpolation,'' 
Lectures at Tianjin University, June 1996.

\bibitem[LS]{ls:in} A. Lascoux and M.-P. Sch\"utzenberger, 
Interpolation de Newton \`a plusieurs variables, in
``S\'eminaire d'Alg\`ebre P.~Dubreil et M.-P.~Malliavin,
36\`eme annie (Paris, 1983--1984),'' Lecture Notes in
Math., Vol.\ 1146, Springer-Verlag, New York, NY, 1985, 161--175.

\bibitem[LR]{lr:gca} D. E. Littlewood and A. R. Richardson,
Group characters and algebra, {\it Philos. Trans. Roy. Soc. London
Ser. A\ } {\bf 233} (1934), 49--141.

\bibitem[M1]{mac:sfh} I. G. Macdonald,  ``Symmetric  functions 
and Hall polynomials,'' 2nd edition, \oup, Oxford, 1995.

\bibitem[M2]{mac:sft} I. G. Macdonald,  Schur functions: theme and
variations,  in ``Actes 28-e S\'eminaire Lotharingien,''
I.R.M.A., Strasbourg, 1992, 5--39.

\bibitem[M3]{mol:fss} A. Molev, Factorial supersymmetric Schur functions
and
super Capelli identities, in ``A.~A.~Kirillov Seminar on
Representation Theory,'' S.~Gindikin, ed., Amer.\ Math.\ Soc.\
Transl.,  Amer.\ Math.\ Soc., Providence, 1997, to appear; q-alg/9606008.

\bibitem[N]{naz:yci} M. Nazarov, Yangians and Capelli identities,
in ``A.~A.~Kirillov Seminar on Representation Theory,''
S.~Gindikin, ed., Amer.\ Math.\ Soc.\ Transl., Amer.\ Math.\
Soc., Providence, 1997, to appear; q-alg/9601027.

\bibitem[O1]{oko:qih} A. Okounkov, Quantum immanants and higher
Capelli identities, {\it Transformation Groups} {\bf 1} (1996), 99--126.

\bibitem[O2]{oko:ybw} A. Okounkov, Young basis, Wick formula, and
higher Capelli identities, {\it Int.\ Math.\ Research Notes} (1996),
817--839.

\bibitem[OO]{oo:ssf} A. Okounkov and G. Olshanski, Shifted Schur functions,
{\it St.\ Petersburg Math.\ J.} {\bf 9} (1997), no. 2; q-alg/9605042.

\bibitem[S1]{sag:sgr} B. E. Sagan, ``The symmetric group: representations,
combinatorial algorithms, and symmetric functions,'' 2nd edition,
Springer-Verlag, New York, to appear.

\bibitem[S2]{sah:sci} S. Sahi, The spectrum of certain invariant
differential
operators associated to a Hermitian symmetric space, in
``Lie Theory and Geometry,'' J.-L. Brylinski, R. Brylinski,
V. Guillemin, V. Kac eds., Progress in Math., Vol.\ 123, Birkh\"auser,
Boston, 1994, 569--576.

\bibitem[V]{vei:t} S. Veigneau, ``Calcul symbolique et calcul
distribu\'e en combinatoire al\-g\'e\-brique,'' Ph.D. thesis,
Universit\'e de Marne-la-Vall\'ee, Marne-la-Vall\'ee, 1996.

\bibitem[Z]{zel:glr} A. V. Zelevinsky, A generalization of the
Littlewood--Richardson rule and the Robinson-Schensted-Knuth 
correspondence, {\it J. Algebra} {\bf 69} (1981), 82--94.

\end{\bib}

\end{document}